\documentclass[a4paper,reqno]{amsart}
\usepackage{geometry}                

\usepackage{a4wide,amssymb}
\usepackage[latin1]{inputenc}
\usepackage[T1]{fontenc}
\usepackage[english]{babel}
\usepackage{amsfonts}
\usepackage{amsmath}
\usepackage{amsthm}
\usepackage{graphicx,xcolor}
\usepackage{latexsym}
\usepackage[colorlinks=true, linkcolor=magenta, citecolor=cyan, urlcolor=blue]{hyperref}
\usepackage{bbm}

\definecolor{light}{gray}{.9}

\def\dl{\frac{d^{\nu}\lambda}{(2\pi)^{\nu}}}
\def\dln{\frac{d^{\nu}\lambda}{(2\pi \sqrt n)^{\nu}}}

\newcommand{\OOO}[1]{O \left(#1\right)}

\newcommand{\oo}[1]{o \left(\frac{1}{#1}\right)}
\newcommand{\ooo}[1]{o \left(#1\right)}

\newcommand{\dd}{\,\text{\rm d}}

\theoremstyle{plain}
\newtheorem{theorem}{Theorem}[section]
\newtheorem{lemma}[theorem]{Lemma}

\theoremstyle{remark}
\newtheorem{remark}[theorem]{Remark}

\def\be{\begin{equation}}
\def\ee{\end{equation}}
\def\bea{\begin{eqnarray}}
\def\eea{\end{eqnarray}}
\def\ni{\noindent}
\def\nn{\nonumber}

\def\diag{\operatorname{diag}}

\def\He{\operatorname{H}}

\def\Si{\operatorname{si}}

\def\a{\alpha}
\def\f{\varphi}

\def\e{\varepsilon}
\def\d{\delta}

\def\r{\rho}

\def\R{\mathbb{R}}

\def\s{\sigma}
\def\a{\alpha}
\def\e{\varepsilon}
\def\d{\delta}

\def\l{\lambda}
\def\r{\rho}
\def\t{\tau}

\def\Z{\mathbb{Z}}

\def\N{\mathbb{N}}

\def\P{\mathbb{P}}
\def\E{\mathbb{E}}

\def\be{\begin{equation}}
\def\ee{\end{equation}}
\def\bc{\begin{center}}
\def\ec{\end{center}}
\def\sign{\operatorname{sign}}
\def\l{\lambda}

\numberwithin{equation}{section}

\DeclareMathSymbol{\leqslant}{\mathalpha}{AMSa}{"36} 
\DeclareMathSymbol{\geqslant}{\mathalpha}{AMSa}{"3E} 
\DeclareMathSymbol{\eset}{\mathalpha}{AMSb}{"3F}     
\renewcommand{\leq}{\;\leqslant\;}                   
\renewcommand{\geq}{\;\geqslant\;}                   

\author{Giuseppe Genovese}
\address{Giuseppe Genovese: Institut f\"ur Mathematik, Universit\"at Z\"urich,
CH-8057 Z\"urich, Switzerland.}
\email{giuseppe.genovese@math.uzh.ch}

\author{Renato Luc\`a}
\address{Renato Luc\`a: Departement Mathematik und Informatik, Universit\"at Basel, Spiegelgasse 1, CH-4051 Basel, Switzerland.}
\email{renato.luca@unibas.ch}

\title{Local Central Limit Theorem for a Random Walk Perturbed in One Point}

\date{\today}

\begin{document}

\maketitle

\begin{abstract}
We consider a symmetric random walk on the $\nu$-dimensional lattice, whose exit probability from the origin is modified by an antisymmetric perturbation and prove the local central limit theorem for this process. A short-range correction to diffusive behaviour appears in any dimension along with a long-range correction in the one-dimensional case.  
\vspace{5mm}

\ni \textbf{MSC:} 60F05, 60G50, 60J10.
\end{abstract}

\section{Statement of the Result}

Let $\{\xi_n\}_{n\in\mathbb N}$ be symmetric i.i.d. random variables taking values 
in $\Z^\nu$, such that 
\be\label{eq:Hyp-xi}
\E[\xi_{1}]=0 \,;
\quad
B_{ij} := \E[\xi_{1,i} \, \xi_{1,j}] \,; 
\quad
\E[\xi_{1,1}^{\alpha_{1}} \dots \, \xi_{1,\nu}^{\alpha_{\nu}}] < \infty, \quad |\alpha|=4\,,
\ee
where $\xi_{1,i}$ is the $i$-th component of $\xi_{1}$ and we use the customary multi-index notation $\alpha:= (\alpha_{1}, \ldots, \alpha_{\nu})\in\N_0^\nu$, $|\alpha| := \alpha_{1} + \ldots + \alpha_{\nu}$. The positive definite matrix $B$ has eigenvalues $\{\s^{2}_{i}\}_{i=1,\dots,\nu}$ and in the one-dimensional case we 
simply write $\sigma$ rather $\sigma_{1}$. We denote by $S_n$ the random walk with increments $\xi_n$ and initial condition $S_0=x_0\in\Z^\nu$ and we require $\{S_n\}_{n\in\mathbb N}$ to be aperiodic and irreducible \cite{law}. 
Moreover let $\{\eta_n\}_{n\in\mathbb N}$ be another sequence of i.i.d. random variables in $\Z^\nu$ such that the matrix $\quad\left(\E[\eta_{1, i} \, \eta_{1, j}]\right)_{i,j=1,\ldots,\nu}$ is positive definite and
\be\label{eq:Hyp-eta}
\E[\eta_{1}]=:d\neq0\,;
\quad \E[\eta_{1,1}^{\alpha_{1}} \dots \, \eta_{1,\nu}^{\alpha_{\nu}}]<\infty, \quad |\alpha|=3\,.
\ee

In this note we consider the aperiodic and irreducible Markov chain $\{X_n\}_{n\in\N}$ on $\Z^\nu$ defined by
\be\label{eq_X_n}
X_n=X_{n-1}+\xi_n+(\eta_n-\xi_n)\d_{X_{n-1},0}\,,\quad X_0=x_0\in\Z^\nu\,,
\ee
where $\d_{x,y}$ denotes as customary the Kronecker delta.
Away from the origin this is just the symmetric random walk. Every time it hits zero, it exits with a different probability given by
$$
\P(X_n=x\,|\,X_{n-1}=0)=\P(\eta_1=x)\,. 
$$
We set
$$
p(x):=\P(\xi_1=x)\,,\quad q(x):=\P(\eta_1=x)\,,
$$
and
$$
a(x):=\P(\eta_1=x)-\P(\xi_1=x)=q(x)-p(x)\,,
$$
so that we can conveniently represent the transition probability for $X_n$ as

\be\label{eq:reprP}
\P(X_n=x\,|\,X_{n-1}=y)=p(x-y)+\d_{y,0}a(x)\,. 
\ee
Of course
$$
\sum_{x\in\Z^\nu} a(x)=0\,,\qquad \sum_{x\in\Z^\nu} xa(x)=d\,.
$$
We will assume antisymmetry of $a(x)$:
\be\label{a-ant}
a(x)=-a(-x)\,.
\ee
Note that this entails
$$
a(x)=q_{a}(x)\,,\qquad p(x)=q_{s}(x)\,,
$$
where $q_{a},q_{s}$ denotes respectively the antisymmetric and symmetric part of $q(x)$. 

The main result of this note follows. As customary, $\oo{n^p}$ denotes a quantity approaching zero faster than $n^{-p}$ uniformly in $x\in\Z^\nu$ and we set $|B|:=\sqrt{\det B}$. 

\begin{theorem}\label{MainThmOLD}
Let $x\in\Z^\nu$, $P_n(x):=\P(X_n=x|X_0=0)$. Then for any $\nu\geq1$ we have $$P_n(0)=\frac{1}{(2\pi n)^\frac{\nu}{2}|B|}+\oo{n^{\frac\nu 2}}\,,$$ while for $x\neq0$ it is
\begin{itemize}
\item
if $\nu=1$
$$
P_n(x)=\frac{e^{-\frac{x^2}{2\s^2n}}}{\sqrt{2\pi\s^2n}}\left(1+\frac{d}{\s^2}\sign (x)\right) 
+\frac{\psi_{1}(x)}{\sqrt n}+\oo{\sqrt n} \,,
$$
where 
$\psi_1$ is an odd function with $\sup_{|x|}  |\psi_1(x)| \leq C$ and $|\psi_1(x)| \leq C \frac{n^{\e}}{|x| }$ for all $\e >0$.
\\
\item If $\nu =2$
$$
P_n(x)=\frac{e^{-\frac{(x,B^{-1}x)}{2n}}}{2\pi|B| n}
+\frac{\psi_{2}(x)}{n}+\oo{n}  \,,
$$
where $\psi_2(x)$ is a bounded odd function such that $|\psi_2(x)| \leq \frac{C}{|x|}$;
\\
\item
if $\nu \geq 3$
$$
P_n(x)=
\frac{e^{-\frac{(x,B^{-1}x)}{2n}}}{(2\pi n)^\frac{\nu}{2}|B|}+\frac{(q_a\ast G_\nu)(x)}{n^{\frac\nu2}}+\oo{n^{\frac\nu 2}}  \, ,
$$
\end{itemize}
where $G_\nu$ is the $\nu$-dimensional Green function.
\end{theorem}

\begin{remark}
The form of the one dimensional short range correction is typical under the assumption that only three moments are finite, see e.g. 
\cite[Theorem 2.3.10]{law}. In higher dimension the short range correction decays faster with $|x|$. Indeed the 
$\nu$-dimensional ($\nu \geq 3$) Green function satisfies $\lim_{|x| \to \infty} |x|^{\nu-2} G_{\nu}(x) \leq C$, see e.g.\cite[Theorem 4.3.1]{law}). 
\end{remark}

Despite the simplicity of this problem, to the best of our knowledge the sole two other works on it are \cite{MZ} and \cite{BP}. In \cite{MZ}  Minlos and Zhizhina proved the local limit theorem for more general perturbation, acting in a finite neighbourhood of the origin, but a.s. bounded increments. 
Boldrighini and Pellegrinotti in \cite{BP} studied the particular case of single point perturbation in one dimension, with analytic increment distribution, giving a more precise information about the terms in the leading order asymptotics. 
Finally for the transient case ($\nu\geq3$) our theorem can be recovered by the more general result \cite[Theorem 7]{siberia}. 

Our approach is simply based on the inversion of the characteristic function, as in the usual local limit theorem for random walks (see e.g. \cite{law}, Chapter 1). We emphasise that since we use real methods, we can drop the analyticity assumptions of \cite{MZ, BP}. Moreover, we do not compute sub-leading corrections to the diffusive scale, but our method in principle would allow this calculation (not without effort, though). In the more involved context of random walk in random environment a similar approach has been used in \cite{C}.

The perturbed walk we consider inherits recurrence or transience from the symmetric random walk. We remark that the antisymmetry assumption on $a(x)$ simplifies the calculations and already captures all the interesting features of the problem. This appears evident for instance looking at the asymptotic formula for $P_n(x)$ in \cite{BP} ($\nu=1$): the symmetric part of the perturbation affects only the form of the short-range correction and the coefficient of the long range correction (not precisely $d/\s^2$, but just proportional to it). 

Lastly, we briefly recall some notions on the Hermite polynomials we will use along the proof (we refer for instance to \cite{bat}). We denote by $\He_\a(x)$ the Hermite polynomials on $\R$. The Rodrigues formula reads 
\be\label{eq:Hermite}
e^{-\frac{x^2}{2}}\He_{\a}(x)=(-1)^{\alpha}\partial^{\alpha}\left(e^{-\frac{x^2}{2}}\right)=\int_{\R}\frac{\dd \l}{\sqrt{2\pi}} 
e^{-ix\l}e^{-\frac{\l^2}{2}}(i\l)^{\alpha}\,.
\ee
The Hermite polynomials form a basis of $L^{2}\left(e^{-\frac{x^{2}}{2}}dx\right)$ with
$$
\int_{\mathbb{R}} \He_{n}(x) \He_{m}(x) \frac{e^{-\frac{x^{2}}{2}}}{\sqrt{2\pi}}dx = n! \delta_{nm}\,.
$$
Any $f\in L^{2}\left(e^{-\frac{x^{2}}{2}}dx\right)$ can be written as
$
f=\sum_{n\geq 0} c_{n} \He_{n}\,,
$
where the equality is in $L^2$ mean and $c_{n}$ are given by
$$
c_{n} = \frac{1}{n!} \int_{\mathbb{R}}   f(x)  \He_{n}(x) \frac{e^{-\frac{x^{2}}{2}}}{\sqrt{2\pi}}dx\, .
$$
In particular, using $\He_n(x)=(-1)^n \He_n(-x)$, for $f=\sign(x)$ we can compute
\be\label{eq:utile}
\frac12\int_{-\infty}^{\infty}\sign(x) \He_{2n+1}(x)\frac{e^{-\frac{x^{2}}{2}}}{\sqrt{2\pi}}dx=\int_0^\infty \He_{2n+1}(x)\frac{e^{-\frac{x^{2}}{2}}}{\sqrt{2\pi}}dx=\He_{2n}(0)=(-1)^n\frac{(2n)!}{2^nn!}\,. 
\ee
The series of $\sign x$ converges pointwise uniformly in each interval \cite[Theorem 9.1.6]{bat}.

Throughout we adopt the following notations: 
$
\widehat{f}(\l):=\sum_{x\in\Z^{\nu}}f(x)e^{i\l x}
$ 
is the usual Fourier transform on $\Z^\nu$; $C$ denotes an absolute constant which may change line by line. 

\subsection*{Acknowledgements} The authors thank C. Boldrighini for suggesting the problem. R. L. is supported by the ERC grant 676675 FLIRT.

\section{Proof}

First we get a nice representation formula for antisymmetric $a$.
\begin{lemma}
Let $a(x)=-a(-x)$ for any $x\in\Z^\nu$. Then
\be\label{eq:rap-antiP}
P_n(x)=p^{\ast n}(x) + \sum_{k=0}^{n-1}p^{\ast k}(0)(q_a\ast p^{\ast (n-k-1)}(x))\,,
\ee
where $p^{\ast n}$ is recursively defined by $p^{\ast n}(x) = (p\ast p^{\ast (n-1)})(x)$,  $p^{\ast 0}(x) = \delta_{x,0}$.
\end{lemma}
\begin{proof}
We introduce the hitting times
$$
\t:=\inf\{n\geq1\,:\,X_n=0\}\,,\quad \t':=\inf\left\{n\geq1\,:\,S_n=0\right\}
$$
and the first return probabilities ($n\geq1$)
$$
f_n(x):=\P(\t= n|X_0=x)\,,\quad f'_n(x):=\P(\t'= n|S_0=x)\,.
$$
We also set
$$
g_n(x,y):=\P(X_n=y,\t\geq n|X_0=x)\,.
$$
Note that for $x\neq 0$ one has
\begin{equation}\label{GnProp}
g_n(x,y):=\P(X_n=y,\t\geq n|X_0=x) = \P(S_n=y,\t\geq n|X_0=x) \, .
\end{equation}
The antisymmetry of $a(y)$ yields $f_1(0)=f'_1(0)=p(0)$ and for $n\geq2$
\bea
f_n(0)&=&\sum_{y\neq0} P_1(y)f_{n-1}(y)=\sum_{y\neq0}\left(p(y)+ a(y)\right)f'_{n-1}(y)\nn\\
&=&\sum_{y \neq 0} p(y)f'_{n-1}(y)+\sum_{y \neq 0} a(y)f'_{n-1}(y)=f'_n(0)\,.\nn
\eea
Therefore we obtain $f_n(0)=f'_n(0)$ for all $n\in\N$, whence $P_n(0)=\P(X_n=0|X_0=0)=\P(S_n=0|S_0=0)=p^{\ast n}(0)$.
Then
\bea
P_n(x)& = & \sum_{k=0}^{n-1}P_k(0)g_{n-k}(0,x)=\sum_{k=0}^{n-1}p^{\ast k}(0)\sum_{y\in\mathbb{Z^{\nu}}/\{0\}}\big(p(y)+ a(y)\big)g_{n-k-1}(y,x)
\nn\\
&=& p^{\ast n}(x)+\sum_{k=0}^{n-1}p^{\ast k}(0)\sum_{y\in\mathbb{Z^{\nu}}/\{0\}}a(y)\P(S_{n-k-1}=x,\t'\geq n-k-1|S_0=y)\nn\,,
\eea
where the first identity is a decomposition with respect to the last visit of the walk to the origin, the second identity follows by the definition of $g_n$
and the last identity follows by \eqref{GnProp}.
Now we use  
\begin{equation}\label{reygfyeg}
\sum_{y\in\mathbb{Z^{\nu}}/\{0\}}a(y)\P(S_{n-k-1}=x,\t'\geq n-k-1|S_0=y)=
\sum_{y\in\mathbb{Z^{\nu}}/\{0\}}a(y)\P(S_{n-k-1}=x|S_0=y)\,,
\end{equation}
which is a consequence of the antisymmetry of $a$.
Indeed we need to show that
\begin{equation}\label{reygfyegreygfyeg}
\sum_{y\in\mathbb{Z^{\nu}}/\{0\}}a(y)\P(S_{n-k-1}=x,\t' < n-k-1|S_0=y) = 0 \, .
\end{equation}
To prove \eqref{reygfyegreygfyeg} we use the symmetry of the walk $S_{n}$ and the antisymmetry of $a$. We have 
\begin{equation}\label{reygfyegreygfyegreygfyeg}
\P(S_{n-k-1}=x,\t' < n-k-1|S_0=y) = \P(S_{n-k-1}=x,\t' < n-k-1|S_0=-y) \, .
\end{equation}
which follows by
\begin{eqnarray}
&& 
\P(S_{n-k-1}=x,\t' < n-k-1|S_0=y) 
 \\ \nonumber
&=&
\sum_{\ell = 1}^{n-k-2} \P(S_{\ell}=0,\t' = \ell|S_0=y) \P(S_{n-k-1}=x|S_{\ell}=0)
\\ \nonumber
&=&
\sum_{\ell = 1}^{n-k-2} \P(S_{\ell}=0,\t' = \ell|S_0=-y) \P(S_{n-k-1}=x|S_{\ell}=0)
\\ \nonumber
&=&
\P(S_{n-k-1}=x,\t' < n-k-1|S_0=-y)  \, ,
\end{eqnarray}
where we used
$$
 \P(S_{\ell}=0,\t' = \ell|S_0=y) = \P(S_{\ell}=0,\t' = \ell|S_0=-y) \, , 
$$
which follows by the symmetry of $\P$.
Combining \eqref{reygfyegreygfyegreygfyeg}
with the antisymmetry of $a$, immediately gives \eqref{reygfyegreygfyeg}.
Thus we get
\bea
P_n(x)&=&p^{\ast n}(x)+\sum_{k=0}^{n-1}p^{\ast k}(0)\sum_{y\in\mathbb{Z^{\nu}}/\{0\}}a(y)\P(S_{n-k-1}=x|S_0=y)\nn\\
&=&\sum_{k=0}^{n-1}p^{\ast k}(0)(a\ast p^{\ast (n-k-1)})(x)\nn\,,
\eea
that is the (\ref{eq:rap-antiP}).
\end{proof}

Let now
\be\label{eq:fourierP}
\f_n(\l):=\sum_{x\in\Z^\nu} e^{i\l x}P_n(x)\,,
\ee
be the characteristic function of $X_n$. Using (\ref{eq:rap-antiP}) it can be expressed as
\bea
\f_n(\l)&=&\widehat{p}^{\,n}(\l)+\phi_n(\l)\,,\label{eq:phi1}\\
\phi_n(\l)&:=&\sum_{k=0}^{n-1}p^{\ast k}(0)\widehat{q_{a}}(\l)\widehat{p}^{\,n-k-1}(\l)\label{eq:phi}\,.
\eea
The formula (\ref{eq:phi}) will be the starting point of our proof. Note that we have already proved that $P_n(0)$ is unchanged by the perturbation, since the second summand on the r.h.s. of the (\ref{eq:rap-antiP}) vanishes in $x=0$. Henceforth we will assume then $x\neq0$.  

Now we set
\begin{eqnarray}\label{IandII}
I_{n}&:=&\sum_{k=0}^{n-1}p^{\ast k}(0)\int_{[-\pi, +\pi]^{\nu}\setminus [-\d_n,\d_n]^{\nu} } \dl e^{-ix \cdot \l} \widehat{p}^{\,n-k-1}(\l) \widehat{q_{a}}(\l),\label{eq:II}\\
II_{n}&:=&\sum_{k=0}^{n-1}p^{\ast k}(0)\int_{[-\d_n,\d_n]^{\nu}} \dl e^{-ix \cdot \l} \widehat{p}^{\,n-k-1}(\l) \widehat{q_{a}}(\l) \label{eq:I}\,,
\end{eqnarray}
where for $\nu=1$ we set 
$\d_n$ a decreasing sequence of positive numbers such that $\lim_{n \to \infty} \delta_{n}^2 \ln n = 0$
and $\d_n:=n^{\e-\frac12}$ with $\e>0$ for $\nu\geq2$. Obviously $P_n(x)=p^{\ast n}(x)+I_n+II_n$.
We evaluate separately $I_n$ and $II_n$ in the following lemmas. 

\begin{lemma}\label{lemma:II} 

\

\begin{enumerate}
\item
Let $\nu=1$. Then 
\begin{equation}\label{eq:I-dim1}
I_{n} = \frac{\psi_1(x)}{\sqrt n}+ \oo{n^{1/2}} \,,
\end{equation}
where 
$\psi_1$ is an odd function that $\sup_{|x|}  |\psi_1(x)| \leq C$ and $|\psi_1(x)| \leq C \frac{n^{\e}}{|x| }$ for all $\e >0$.
\\
\item   
Let $\nu=2$. Then 
\begin{equation}\label{eq:I-dim2}
I_{n}= \frac{\psi_2(x)}{n} + \oo{n} \, ,
\end{equation}
where $\psi_2$ is a bounded odd function such that $\lim_{|x| \to \infty} |x| |\psi_2(x)| \leq C$.
\\
\item
Let $\nu\geq3$ and $G_\nu$ denote the $\nu$-dimensional Green function. Then 
\begin{equation}\label{eq:I-dim3}
I_{n}=\frac{(q_a\ast G_\nu)(x)}{n^{\nu/2}}+\oo{n^{\nu/2}} \, .
\end{equation}
\end{enumerate}
\end{lemma}

\begin{lemma}\label{lemma:I} 

\

\begin{enumerate}
\item
Let $\nu=1$. Then 
\begin{equation}\label{eq:dim1}
II_{n}=\frac{e^{-\frac{x^2}{2\s^2}}}{\sqrt{2\pi\s^2n}}\frac{d}{\s^2}\sign (x)+\oo{n^{1/2}} \,. 
\end{equation}

\

\item
Let $\nu\geq2$. Then 
\begin{equation}\label{eq:dim3}
II_{n}=\oo{n^{\nu/2}} \, .
\end{equation}
\end{enumerate}
\end{lemma}

According to formulas (\ref{eq:fourierP}) and (\ref{eq:phi1}), we can then prove the main theorem by Fourier inversion using the local central limit theorem for the homogeneous probability \cite{law}. 
The proofs of the lemmas are presented below.

\begin{proof}[Proof of Lemma \ref{lemma:II}]

We recall that for $k \geq1$
\be\label{eq:p0} 
p^{\ast k}(0) = \frac{1}{(2\pi|B|k)^{\frac\nu2}}+\oo{k^{\frac\nu2}}\,.
\ee 
and due to aperiodicity there exists $b>0$ such that 
\begin{equation}\label{eq:p1}
|\hat p(\l)|\leq 1-b|\l|^2 \leq e^{-b|\l|^2}\,
\end{equation}
for all $\l \in [-\pi, \pi]^\nu$ (see for instance Lemma 2.3.2 in \cite{law}).
We fix $\theta \in (0,1)$ and we split the sum over $0 \leq k \leq n-n^{\theta}$ and $n-n^{\theta} < k \leq n-1$. By \eqref{eq:p0} and \eqref{eq:p1} we get
\be
\left| \sum_{0 \leq k \leq n-n^{\theta}} p^{\ast k}(0)\hat p^{n-k-1}(\l) \right| 
\leq
 \sum_{n^{\theta} -1  \leq h \leq n-1} p^{\ast (n-h-1)}(0) | \hat p(\l)| ^{h} 
\leq e^{-b n^{\theta}\l^2}\r_{n}\,,
\ee
where $\r_{n}:=\sum_{n^{\theta} \leq h \leq n}p^{\ast h}(0)$ that, recalling \eqref{eq:p0}, is 
bounded by a constant (uniform in $n$) for $\nu\geq3$, by $C \ln n $ for $\nu = 2$ and by $C \sqrt{n}$ for $\nu =1$. 
Again by \eqref{eq:p0} we get 
\bea
\sum_{n-n^{\theta} < k \leq n-1} p^{\ast k}(0)\hat p^{n-k-1}(\l)&=&\sum_{0 \leq h  <  n^{\theta} - 1} p^{\ast (n-h-1)}(0)\hat p(\l)^{h}\nn\\
& = &\frac{1 + o(1)}{|B|(2\pi n)^{\nu/2}}\sum_{h=0}^{h= \lfloor n^{\theta} - 1 \rfloor} \hat p(\l)^{h}\nn\\
&=&\frac{1 + o(1)}{|B|(2\pi n)^{\nu/2}}\frac{1-\hat p^{ \lfloor n^{\theta}  \rfloor}(\l)}{1-\hat p(\l)}\,.
\eea
Therefore
$$
I_n = 
R_1(\l)  + \frac{1 + o(1)}{|B|(2\pi n)^{\nu/2}}\int_{[-\pi, +\pi]^{\nu}\setminus [-\d_n,\d_n]^{\nu}} \dl e^{-i x \cdot \l } 
\frac{1-\hat p^{\lfloor n^{\theta} \rfloor}(\l)}{1-\hat p(\l)}\widehat{q_{a}}(\l)\,
$$
where 
$$
|R_1(\l) | \leq e^{- b n^{\theta} \d_n^2}\r_{n}\sup_{|\l|>\d_n}\widehat{q_{a}}(\l) \leq C e^{- b n^{\theta} \d_n^2}\r_{n} \, .
$$
We can further decompose 
\be
 \int_{[-\pi, +\pi]^{\nu}\setminus [-\d_n,\d_n]^{\nu}} \dl e^{-i x \cdot \l } 
\frac{1-\hat p^{n^{\theta} }(\l)}{1-\hat p(\l)}\widehat{q_{a}}(\l) 
=  \int_{[-\pi, +\pi]^{\nu}\setminus [-\d_n,\d_n]^{\nu}} \dl e^{-i x \cdot \l } 
\frac{\widehat{q_{a}}(\l)}{1-\hat p(\l)} 
+ R_2(\l) \,. 
\ee
where
$$
|R_2(\l)| \leq e^{- b n^{\theta}\d_n^2}\left(\max_{|\l|>\d_n}\frac{\widehat{q_{a}}(\l)}{1-\hat p(\l)}\right) 
\leq e^{- b n^{\theta}\d_n^2} \d_n^{-2}\left(\max_{|\l|>\d_n}\widehat{q_{a}}(\l)\right) 
\leq C e^{- b n^{\theta}\d_n^2} \d_n^{-2} 
\,.
$$
Summarising 
\be\label{eq:I-main}
I_n=\frac{1}{|B|(2\pi n)^{\nu/2}}\int_{[-\pi, +\pi]^{\nu}\setminus [-\d_n,\d_n]^{\nu}} \dl e^{-i x \cdot \l } 
\frac{\widehat{q_{a}}(\l)}{1-\hat p(\l)}+\OOO{ e^{-n^{\theta}\d_n^2}(\d_n^{-2}+\r_{n})} + \oo{n^{\nu/2}}\,.
\ee
$I_n$ is an odd function of $x\in\Z^{\nu}$ (due to the parity of $\hat p, \hat q_a$) plus a remainder $\oo{n^{\nu/2}}$.
We complete the proof treating separately the transient and the recurrent case.

\begin{proof}[Proof of \eqref{eq:I-dim3} ($\nu\geq3$)] Let first $\nu\geq3$ and recall that 
$$
\int_{[-\pi, +\pi]^{\nu}} \dl 
\frac{e^{-i x \cdot \l } }{1-\hat p(\l)} = G_\nu(x) 
$$
is the Green function of the free walk. We can write
\bea\nonumber
\int_{[-\pi, +\pi]^{\nu}\setminus [-\d_n,\d_n]^{\nu}} \dl e^{-i x \cdot \l } 
\frac{\widehat{q_{a}}(\l)}{1-\hat p(\l)}&=&
(q_a\ast G_\nu)(x) -
\int_{[-\d_n,\d_n]^{\nu}} \dl e^{-i x \cdot \l }  \frac{\widehat{q_{a}}(\l)}{1-\hat p(\l)} \, .
\eea
Now we show that  
\begin{equation}\label{RefEq}
\int_{[-\d_n,\d_n]^{\nu}} \dl e^{-i x \cdot \l }  \frac{\widehat{q_{a}}(\l)}{1-\hat p(\l)} = \OOO { \delta_n^{\nu - 1} },
\quad 
\nu \geq 2 \, ,
\end{equation}
so that the statement (3) follows 
letting $\theta = 1/2$ since $\delta_{n} = n^{\e - 1/2 }$ for some $\e >0$. 
To prove \eqref{RefEq} we use \ref{eq:p1} and $|\widehat{q_{a}}(\l)| \leq C |\lambda|$ for $\delta_n$ (and so $\l$) sufficiently small.
This leads us to 
\begin{eqnarray}
\left| \int_{[-\d_n,\d_n]^{\nu}} \dl e^{-i x \cdot \l }  \frac{\widehat{q_{a}}(\l)}{1-\hat p(\l)} \right|
& \leq & C_{\nu}  \int_{[-\d_n,\d_n]^{\nu}} \dl   \frac{1}{|\l|} \leq C_{\nu} \delta_n^{\nu - 1}, \quad \nu \geq 2 \, . 
\\ \nonumber
& &
\end{eqnarray}
\end{proof}
\

In the recurrent cases we have to evaluate directly the integral in (\ref{eq:I-main}). 

\begin{proof}[Proof of (\ref{eq:I-dim2}) ($\nu=2$)]
Let us assume $x_1 \geq x_2$ (the case $x_2 \geq x_1$ is analogous). Noting that 
$\frac{\widehat{q_{a}}(\l)}{1-\hat p(\l)}$ is everywhere regular but in the origin,
we decompose
\begin{eqnarray}\label{THTFTOTRHS}
 \int_{[-\pi, +\pi]^{2}\setminus[-\delta_n, \delta_n]^2}  \frac{d^{2}\l}{ (2 \pi)^2} e^{ - i x \cdot \l } 
\frac{\widehat{q_{a}}(\l)}{1-\hat p(\l)}
& = &
\int_{[-\pi, +\pi]^{2}} \frac{d^{2}\l}{ (2 \pi)^2} e^{ - i x \cdot \l }
\frac{2id \cdot \l}{(B \l, \l)}  
\\  \nonumber
&+ & 
\int_{[-\pi, +\pi]^{2}} \frac{d^{2}\l}{ (2 \pi)^2}  e^{ - i x \cdot \l } 
 \left(  \frac{\widehat{q_{a}}(\l)}{1-\hat p(\l)} - \frac{2id \cdot \l}{(B \l, \l)}  \right)  +  \OOO{\delta_n} \, ,
\end{eqnarray}
where we used again \eqref{RefEq}.
Now we show that 
$\partial_{\l_1} \left( \frac{\widehat{q_{a}}(\l)}{1-\hat p(\l)} - \frac{2id \cdot \l}{(B \l, \l)} \right)$ is continuous on $[-\pi, \pi]^2$. 
Recalling the parity of $\widehat{q_{a}}(\l)$
and $\hat p(\l)$ we can expand 
$$
\widehat{q_{a}}(\l) = id\cdot \lambda + r_1(\l), 
\qquad
\hat p(\l) = 1 - \frac12 (B \l, \l) + r_2(\l)\,,
$$
where $r_1(\l) = \ooo{\lambda^3}$ and $r_2(\l) = \ooo{\lambda^4}$. Then
we can write $\partial_{\l_1} \left( \frac{\widehat{q_{a}}(\l)}{1-\hat p(\l)} - \frac{2id \cdot \l}{(B \l, \l)} \right)|_{\lambda =0}$ as
$$
 \lim_{\l\to0}\frac{1}{\l_1}\left(\frac{\widehat{q_{a}}(\l)}{1-\hat p(\l)} - \frac{2id \cdot \l}{(B \l, \l)}\right)= \lim_{\l\to0}\frac{1}{\l_1}\left(\frac{2r_1(\l)(B\l,\l)+4id\cdot\l r_2(\l))}{(B\l,\l)((B\l,\l)-2r_2(\l))}\right)\,,
$$
that is bounded. Thus, integrating by parts 
the function $e^{ - i x \cdot \l } =  \frac{i}{ x_1} \partial_{\l_1} e^{ - i x \cdot \l}$ in the variable $\l_1$, we get
\begin{equation}
\left| \int_{[-\pi, +\pi]^{2}} \frac{d^{2}\l}{ (2 \pi)^2}  e^{ - i x \cdot \l }  \left( 
\frac{\widehat{q_{a}}(\l)}{1-\hat p(\l)} - \frac{2i d \cdot \l}{(B \l, \l)}  \right)  \right| \leq \frac{C}{|x_1|} \leq \frac{C \sqrt{2}}{|x|} \,.
\end{equation}
To handle the first term on the r.h.s. of \eqref{THTFTOTRHS} we first notice that, again integrating by parts
and using the continuity of
$\frac{d \cdot \l}{(B \l, \l)} $ on $[-\pi, \pi]^2\setminus\{|\l|\leq \pi\}$, we get
\begin{equation}\label{THTFROTRHSOT}
\int_{[-\pi, +\pi]^{2}} \frac{d^{2}\l}{ (2 \pi)^2} e^{ - i x \cdot \l }
\frac{d \cdot \l}{(B \l, \l)}
=
\frac{1}{n} 
\int_{|\lambda| \leq \pi} \frac{d^{2}\l}{ (2 \pi)^2} e^{-i x \cdot \l } 
\frac{d \cdot \l}{(B \l, \l)} +  R_1(x)
\end{equation}
with
$
|R_1(x)| \leq \frac{C}{|x|} \, 
$.
We change variables 
$\l = O \bar \l$ and $x= \bar x O^{-1}$, with $O^{-1} = O^{T}$ and $O^{-1} B O=\diag(\sigma_{1}^2,\s_2^2)$, so that
\be
\int_{|\l| \leq \pi} \frac{d^{2}\l}{ (2 \pi)^2} e^{- i x \cdot \l}  \frac{d \cdot \l}{(B \l, \l)}=
\int_{|\bar \l| \leq \pi} \frac{d^{2} \bar \l}{ (2 \pi)^2} e^{- i \bar x \cdot \bar \l  }
\frac{\bar d_1 \bar \l_1 + \bar d_2 \bar \l_2  }{\sigma_1^2 \bar \l_1^2 + \sigma_2^2 \bar \l_2^2}  \, ,
\ee
where $\bar d = O d$.
Since the function $\frac{\bar d_1 \bar \l_1 + \bar d_2 \bar \l_2  }{\sigma_1^2 \bar \l_1^2 + \sigma_2^2 \bar \l_2^2}$
is continuous on $\{|\l \leq \pi|\} \setminus [-\pi/\sqrt{2}, \pi / \sqrt{2} ]^2$, integrating by parts 
we get  
\begin{equation}\label{FinalmQST}
\int_{|\l| \leq \pi} \frac{d^{2}\l}{ (2 \pi)^2} e^{-i x \cdot \l}  
\frac{d \cdot \l}{(B \l, \l)}
= 
\int_{[-\pi / \sqrt{2}, \pi / \sqrt{2}]^{2} } \frac{d^{2} \bar \l}{ (2 \pi)^2} e^{-i \bar x \cdot \bar \l  }
\frac{\bar d_1 \bar \l_1 + \bar d_2 \bar \l_2  }{\sigma_1^2 \bar \l_1^2 + \sigma_2^2 \bar \l_2^2} 
+ R_2(x) \, ,
\end{equation}
with 
$
|R_2(x)| \leq\frac{C }{|x|} \,
$.
Then we will prove that also the first term on the r.h.s. of \eqref{FinalmQST} satisfies the same bound
$$
\left| \int_{[-\pi / \sqrt{2}, \pi / \sqrt{2}]^{2} } \frac{d^{2} \bar \l}{ (2 \pi)^2} e^{-i \bar x \cdot \bar \l  }
\frac{\bar d_1 \bar \l_1 + \bar d_2 \bar \l_2  }{\sigma_1^2 \bar \l_1^2 + \sigma_2^2 \bar \l_2^2} \right| \leq \frac{C}{|x|} \, .
$$
By triangle inequality and symmetry we only need to prove
$$
\left| \int_{[-\pi / \sqrt{2}, \pi / \sqrt{2}]^{2} } \frac{d^{2} \bar \l}{ (2 \pi)^2} e^{-i \bar x \cdot \bar \l  }
\frac{\bar d_1 \bar \l_1  }{\sigma_1^2 \bar \l_1^2 + \sigma_2^2 \bar \l_2^2} \right| \leq \frac{C}{|x|} \, .
$$
We use
$$
\frac{1}{\s_1^2 \bar \l_1^2+\s_2^2 \bar \l_2^2}=\int_0^\infty dye^{-y(\s_1^2 \bar \l_1^2+  \s_2^2 \bar \l_2^2)} 
$$
and Fubini's theorem to obtain
\begin{eqnarray}\nonumber
& & 
\int_{[-\pi / \sqrt{2}, \pi / \sqrt{2}]^{2} } \frac{d^{2} \bar \l}{ (2 \pi)^2} e^{-i \bar x \cdot \bar \l  }
\frac{\bar d_1 \bar \l_1   }{\sigma_1^2 \bar \l_1^2 + \sigma_2^2 \bar \l_2^2} 
\\ \nonumber
& = &
 \bar d_1\int_0^\infty dy \left[ \int_{-\pi}^\pi  \frac{d \bar \l_2}{2\pi} e^{-i \bar x_2 \bar \l_2- y \s_2^2 \bar \l_2^2}\right]
 \left[\int_{-\pi}^\pi \frac{d \bar \l_1}{2\pi} 
 e^{-i \bar x_1 \bar \l_1- y\s_1^2 \bar \l_1^2} \bar \l_1\right]
\\ \nonumber
&=& 
 \bar d_1 
 \int_0^\infty dy \left[ \int_{\mathbb{R}}  \frac{d \bar \l_2}{2\pi} e^{-i \bar x_2 \bar \l_2- y \s_2^2 \bar \l_2^2}\right]
 \left[\int_{\mathbb{R}} \frac{d \bar \l_1}{2\pi} 
 e^{-i \bar x_1 \bar \l_1- y\s_1^2 \bar \l_1^2} \bar \l_1\right]
\\ \nonumber
&+ & 
 \bar d_1\int_0^\infty dy \left[ \int_{ (-\infty, -\pi) \cup (\pi, \infty) }  \frac{d \bar \l_2}{2\pi} e^{-i \bar x_2 \bar \l_2- y \s_2^2 \bar \l_2^2}\right]
 \left[\int_{(-\infty, -\pi) \cup (\pi, \infty)} \frac{d \bar \l_1}{2\pi} 
 e^{-i \bar x_1 \bar \l_1- y\s_1^2 \bar \l_1^2} \bar \l_1\right]
\\ \nonumber
&-& 
 \bar d_1 
 \int_0^\infty dy \left[ \int_{\mathbb{R}}  \frac{d \bar \l_2}{2\pi} e^{-i \bar x_2 \bar \l_2- y \s_2^2 \bar \l_2^2}\right]
 \left[\int_{(-\infty, -\pi) \cup (\pi, \infty)} \frac{d \bar \l_1}{2\pi} 
 e^{-i \bar x_1 \bar \l_1- y\s_1^2 \bar \l_1^2} \bar \l_1\right]
\\ \nonumber
&-& 
 \bar d_1 
 \int_0^\infty dy \left[ \int_{ (-\infty, -\pi) \cup (\pi, \infty) }  \frac{d \bar \l_2}{2\pi} e^{-i \bar x_2 \bar \l_2- y \s_2^2 \bar \l_2^2}\right]
 \left[\int_{\mathbb{R}} \frac{d \bar \l_1}{2\pi} 
 e^{-i \bar x_1 \bar \l_1- y\s_1^2 \bar \l_1^2} \bar \l_1\right]
\end{eqnarray}
The first term can be explicitly computed and its modulus equals
\begin{equation}\label{TMOTO}
\left|
\frac{1}{4 \s_1 \s_2} \left( \frac{  d_1 \bar x_1}{\s_1^2}  \right)\int_0^{\infty} \frac{dy}{y^2} e^{- \frac{1}{4 y}\left(\frac{ \bar x_1^2}{\s_1^2} + \frac{ \bar x_2^2}{\s_2^2}\right) } 
\right| 
\leq \frac{C}{|x|}\,,
\end{equation}
where the last inequality follows by the change of variables $y' = \left(\frac{ \bar x_1^2}{\s_1^2} + \frac{ \bar x_2^2}{\s_2^2}\right)^{-1} y$. 
For the second term we will show that
\begin{equation}\label{SSDEC}
\left| 
\int_0^\infty dy \left[ \int_{\pi}^{\infty}  d \bar \l_2 e^{-i \bar x_2 \bar \l_2- y \s_2^2 \bar \l_2^2}\right]
 \left[\int_{\pi}^{\infty}  d \bar \l_1 
 e^{-i \bar x_1 \bar \l_1- y\s_1^2 \bar \l_1^2} \bar \l_1\right] \right| \leq \frac{C}{|x_1| |x_2|} \leq \frac{2 C}{|x|} \,.
\end{equation}
Then, since the other contributions are analogous, the statement will follows by triangle inequality.
Integrating by parts the functions $e^{-i \bar \lambda_j \bar x_j } = \frac{i}{\bar x_j} \partial_{\bar \l_j} e^{-i \bar \lambda_j} \bar x_j$ we can bound 
\be\label{TPOT1}
\left| \int_{ \pi}^{\infty}  d \bar \l_2 e^{-i \bar x_2 \bar \l_2- y \s_2^2 \bar \l_2^2} \right|
\leq
\frac{1}{|\bar x_2|} \left( e^{- y \s_2^2 \pi^2} + \int_{\pi}^{\infty}  d \bar \l_2 2 y \s_2^2 \bar \l_2 e^{- y \s_2^2 \bar \l_2^2}  \right)
\leq \frac{2}{|x_2|} e^{- y \s_2^2 \pi^2} 
\ee
and similarly
\begin{eqnarray}\nonumber
\left| \int_{ \pi}^{\infty}  d \bar \l_1 e^{-i \bar x_1 \bar \l_1- y \s_1^2 \bar \l_1^2} \bar \l_1 \right|
& \leq & 
\frac{1}{|\bar x_1|} \left( 
 e^{- y \s_1^2 \pi^2} 
 +
\int_{\pi}^{\infty} d \bar \l_1  2 y \s_1^2 \bar \l_1^2 e^{- y \s_1^2 \bar \l_1^2} 
+
\int_{\pi}^{\infty}  d \bar \l_1 e^{- y \s_1^2 \bar \l_1^2}
 \right)
\\ \nonumber
& \leq &
\frac{1}{|\bar x_1|} \left( e^{- y \s_1^2 \pi^2} 
+  
\frac{2}{\sqrt{y} \s_1} \int_{\sqrt{y} \s_1 \pi}^{\infty} d \bar \eta_1    \bar \eta_1^2 e^{- \bar \eta_1^2} 
+
\frac{2}{\sqrt{y} \s_1}
\int_{\sqrt{y} \s_1 \pi}^{\infty}  d \bar \eta_1 e^{-  \bar \eta_1^2} \right)\nn\\
&\leq &
\frac{1}{| x_1|} \left( e^{- y \s_1^2 \pi^2}  + \frac{C e^{- y \frac{\pi^2 \sigma_1^2}{2}}}{\sqrt{y} \s_1} \right)\label{TPOT2}\,,
\end{eqnarray}
where in the last lines we changed variables $\bar \eta_1 = \sqrt{y} \s_1 \bar \l_1$ and used 
$\int_{0}^{\infty} d \bar \eta_1    \bar \eta_1^2 e^{- \bar \eta_1^2} \leq C/2$ and  
$\int_{0}^{\infty} d \bar \eta_1    \bar  e^{- \bar \eta_1^2} \leq C/2$. Since the product of 
\eqref{TPOT1} and \eqref{TPOT2} is integrable over $y \in (0,\infty)$ the \eqref{SSDEC} has been proved.
Turning to the third term, we notice that it equals
\begin{equation}\nonumber
\left|
\frac{1}{2 \sqrt{\pi} \s_2} \int_0^{\infty} \frac{dy}{\sqrt{y}} e^{- \frac{\bar x_2^2}{4 y \sigma_2^2}   } 
\left[\int_{(-\infty, -\pi) \cup (\pi, \infty)} \frac{d \bar \l_1}{2\pi} 
 e^{-i \bar x_1 \bar \l_1- y\s_1^2 \bar \l_1^2} \bar \l_1\right]
\right|\,. 
\end{equation}
Recalling \eqref{TPOT2}, this can be bounded by
\begin{equation}\nonumber
\frac{C}{| x_1|} 
\int_0^{\infty} \frac{dy}{\sqrt{y}} e^{- \frac{\bar x_2^2}{4 y \sigma_2^2}   } 
 \left( e^{- y \s_1^2 \pi^2}  + \frac{C e^{- y \frac{\pi^2 \sigma_1^2}{2}}}{\sqrt{y} \s_1} \right)
 \leq \frac{C}{| x_1|} \, .
\end{equation}
The estimate for the fourth term is analogous.
\end{proof}

\

\begin{proof}[Proof of \eqref{eq:I-dim1} ($\nu=1$)]

The leading term in (\ref{eq:I-main}) reads
$$
\frac{1}{\sqrt{n}}\int_{[-\pi,\pi]\setminus[-\d_n,\d_n]} \frac{d\l}{2\pi} e^{ix\l}
\frac{\widehat{q_{a}}(\l)}{1-\hat p(\l)}=\frac{1}{\sqrt{n}}\int_{[-\pi,\pi]\setminus[-\d_n,\d_n]} \frac{d\l}{2\pi} \left(\frac{de^{ix\l}}{\s^2\l}+
e^{ix\l} \left( \frac{\widehat{q_{a}}(\l)}{1-\hat p(\l)}-\frac{id}{\s^2\l}\right)\right)\,.
$$
The derivative of the integrand function in the second summand on the r.h.s. is bounded on $[-\pi,\pi]$ (the proof is as in the $\nu=2$ case), therefore the integral can be estimated by $C/|x|$. The first integral gives the main contribution:
\be
\frac{d}{\s^2\pi\sqrt{ n}}\int_{\d_n}^\pi d\l \frac{\sin(x  \l) }{\l}=\frac{d}{\s^2\pi\sqrt{ n}}\left(\Si(\d_n x)-\Si(\pi x)\right)\,,
\ee
where $\Si(t) = \int_{t}^{\infty} \frac{\sin(s)}{s} ds$ is the sine integral function. It stays bounded for small argument and decays as $1/t$ for large $t$. This gives \eqref{eq:I-dim1}. 
\end{proof}
This concludes the proof of the first lemma. 
\end{proof}


\begin{proof}[Proof of Lemma \ref{lemma:I}]

The change of variables $\l\mapsto \l\sqrt{n}$ maps the cube $[-\d_n,\d_n]^\nu$ in $[-\d_n\sqrt n,\d_n\sqrt n]^\nu$. By the cumulant expansion of $\widehat{p}(\l)$ and the Taylor expansion of $\widehat{q_{a}}(\l)$ we obtain
\begin{equation}\label{IIn-long}
II_n=\int_{[-\d_n \sqrt{n},\d_n \sqrt{n}]^{\nu}} \dln e^{-i\frac{x}{\sqrt n}\l}e^{-\frac{(B\l,\l)}{2}}\sum_{k=0}^{n-1}p^{\ast k}(0)e^{\frac{(B\l,\l)k}{2n}}\left[\frac{i(\l\cdot d)}{\sqrt n}+r\left(\frac{\l}{\sqrt n}\right)\right]\,,
\end{equation}
with $\limsup_{|\l|\to0}|r(\l)|/|\l|^3 \leq C$. 
Now we separately handle the case $\nu \geq 2$, and $\nu =1$.

\begin{proof}[Proof of (\ref{eq:dim3}) ($\nu \geq 2$)]
Because of our choice of $\d_n$ for $\nu\geq 2$ we have
\bea
|II_n|&=&\left|\int_{[-n^\e,n^\e]^{\nu}} \dln e^{-i\frac{x}{\sqrt n}\l}e^{-\frac{(B\l,\l)}{2}}\sum_{k=0}^{n-1}p^{\ast k}(0)e^{\frac{(B\l,\l)k}{2n}}\left[\frac{i(\l\cdot d)}{\sqrt n}+r\left(\frac{\l}{\sqrt n}\right)\right]\right|\,\nn\\
&\leq & C \int_{[-n^\e,n^\e]^{\nu}}  \dln \left[\frac{i(\l\cdot d)}{\sqrt n}+r\left(\frac{\l}{\sqrt n}\right)\right]\leq C' \frac{n^{\e(\nu+1)}}{\sqrt{n}n^{\frac\nu2}}=\oo{n^{\frac\nu2}}\,.
\eea
\end{proof}


\begin{proof}[Proof of (\ref{eq:dim1}) ($\nu = 1$)]

It is convenient to set for $|\l|\leq \d_n\sqrt n$
$$
K_{n}(\l):= \frac{\sqrt{n}}{2 \s^2\l^2 \sqrt{\pi}}   \int_{\frac{\s^2\l^2}{2n}}^{\frac{\s^2\l^2(n-1)}{2n}}dt \frac{e^t}{\sqrt t}\,.
$$
We can write $K_{n}(\l)$ as a series as follows
\begin{equation}
K_{n}(\l)=\sum_{\ell\geq0} \left(\frac{\s^2\l^2}{2n}\right)^{\ell}\frac{1}{\ell!}\frac{(n-1)^{\ell+\frac12}-1}{\ell+\frac12}\label{eq:imporant}\,.
\end{equation}
The proof is done simply by expanding the exponential. Then
$$
\sum_{k=1}^{n-1}p^{\ast k}(0) e^{\frac{\s^2\l^2k}{2n}}=\frac{1}{\sqrt{2 \pi \s^2n}}\sum_{t=1/n}^{(n-1)/n}\frac{e^{\frac{\s^2\l^2t}{2}}}{\sqrt{t}}+\oo{\sqrt{n}}\,
$$
and estimating the Riemann sum with the integral (through the second derivative) we obtain
\be
\left|\frac{1}{\sqrt {2\pi\s^2n}}\sum_{t=1/n}^{(n-1)/n}\frac{e^{\frac{\s^2\l^2t}{2}}}{\sqrt{t}}-K_{n}(\l)\right|
\leq\frac{e^{\frac{\s^2\l^2}{2}\left(1-\frac1n\right)}}{n^{\frac32}}+\frac{e^{\frac{\s^2\l^2}{2n}}}{n}\,.
\ee
When we plug this last expression into (\ref{IIn-long}), this can be readily estimated by
$$
\frac{1}{n^{\frac32}}\int_{\R} \frac{d\l}{2\pi \sqrt{n}}e^{-\frac{\s^2\l^2}{2n}}\left[\frac{i(\l  d)}{\sqrt n}+r\left(\frac{\l}{\sqrt n}\right)\right]+\frac{d}{n}\int_{\R} \frac{d\l}{\sqrt{2\pi n}}e^{-\frac{\s^2\l^2}{2}\left(1-\frac1n\right)}\left[\frac{|\l|}{\sqrt n}+r\left(\frac{\l}{\sqrt n}\right)\right]
=\oo{\sqrt n}\,.\nn
$$
Therefore
\begin{equation}\label{WePutHereSp}
II_n=\int_{[-\d_n \sqrt{n},\d_n \sqrt{n}]} \frac{d\l}{2\pi \sqrt{ n}} e^{-i\frac{x}{\sqrt n}\l}e^{-\frac{\s^2\l^2}{2}}K_{n}(\l)\left[\frac{i(\l  d)}{\sqrt n}+r\left(\frac{\l}{\sqrt n}\right)\right]+\oo{\sqrt n}\,.
\end{equation}
Using $|r(\l)|\leq C |\l|^3$ for $\l\in[-\d_n,\d_n]$ we have
\bea
&&\int_{[-\d_n \sqrt{n},\d_n \sqrt{n}]} \frac{d\l}{2\pi \sqrt{ n}} e^{-i\frac{x}{\sqrt n}\l}e^{-\frac{\s^2\l^2}{2}}K_{n}(\l)r\left(\frac{\l}{\sqrt n}\right)\nn\\
&\leq& \int_{[-\d_n \sqrt{n},\d_n \sqrt{n}]} \frac{d\l}{2\pi \sqrt{ n}}e^{-\frac{\s^2\l^2}{2}}\frac{K_{n}(\l)|\l|}{\sqrt n}\frac{|\l|^2}{n}\nn\\
&=& 2\int_{[-\d_n \sqrt{n},\d_n \sqrt{n}]} \frac{d\l}{2\pi \sqrt{ n}}e^{-\frac{\s^2\l^2}{2}}\sum_{\ell\geq0}\frac{(\s|\l|)^{2\ell}}{\ell!(2\ell+1)2^\ell}\frac{|\l|^3}{n}\,.
\eea
Now we split the sum with $n$; the tail is easily bounded:
\be
2\int_{[-\d_n \sqrt{n},\d_n \sqrt{n}]} \frac{d\l}{2\pi \sqrt{ n}}e^{-\frac{\s^2\l^2}{2}}\sum_{\ell\geq n}\frac{(\s|\l|)^{2\ell+1}}{\ell!(2\ell+1)2^\ell}\frac{\l^2}{n}\nn\\
\leq\frac{2}{2n+1}\int_{[-\d_n \sqrt{n},\d_n \sqrt{n}]} \frac{d\l}{2\pi \sqrt{ n}}\frac{|\l|^3}{n}=\oo{\sqrt{n}}\,.\nn
\ee
For the remaining part of the sum we have
\bea
2\int_{[-\d_n \sqrt{n},\d_n \sqrt{n}]} \frac{d\l}{2\pi \sqrt{ n}}e^{-\frac{\s^2\l^2}{2}}\sum_{\ell=0}^n\frac{(\s|\l|)^{2\ell+1}}{\ell!(2\ell+1)2^\ell}\frac{\l^2}{n}
&\leq&\sum_{\ell=0}^n\frac{2\d^2_n}{\sqrt n\ell!(2\ell+1)2^\ell}\int_{\R} \frac{d\l}{2\pi}e^{-\frac{\s^2\l^2}{2}}(\s|\l|)^{2\ell+1}\nn\\
&=&\sum_{\ell=0}^n\frac{1}{2\ell+1}\frac{2\d^2_n}{\sqrt n}=\frac{2\d^2_n\ln n}{\sqrt n}=\oo{\sqrt n}\,,\nn
\eea
due to our choice of $\d_n$ such that $\limsup_{n \to \infty} \delta_{n}^2 \ln n = 0$.
For the linear term in $\l$ we get
\bea
II_n&=&2\sum_{\ell\geq0}\int_{-\d_n\sqrt n}^{\d_n\sqrt n}  \frac{d\l}{ 2 \pi \sqrt{ n}}e^{-i\frac{x}{\sqrt n}\l}e^{-\frac{\l^2\s^2}{2}
} 
\left(    \left(\frac{\l^2\s^2}{2 }\right)^{\ell} \frac{1}{\ell! \, (2\ell+1)}  \right)
\frac{i\l d}{\sqrt n}+\oo{\sqrt n}\label{WePutHereSp}\nn\\
&=&\frac{d}{\s^2}2\sum_{\ell\geq0}\frac{1}{2^\ell\ell! \,(2\ell+1)}\int_{\R}  \frac{d\l}{ 2\pi \sqrt{n} }e^{-i\frac{x}{\s\sqrt n}\l}e^{-\frac{\l^2}{2}
}     \l^{2\ell} (i\l)+\oo{\sqrt n}\nn\\
&=&\frac{d}{\s^2}2\sum_{\ell\geq0}\frac{(-1)^\ell}{2^\ell\ell! \,(2\ell+1)}\frac{e^{-\frac{x^2}{2\s^2}}}{2\pi \sqrt{ n\s^2}}H_{2\ell+1}\left(\frac{x}{\s\sqrt n}\right)+\oo{\sqrt n}\nn\,.
\eea
The leading term can be explicitly computed. Using (\ref{eq:utile}) we obtain the expansion in Hermite polynomials of $\sign \left( \frac{x}{\s\sqrt n} \right)=\sign (x)$: 
\bea
II_n&=&\frac{d}{\s^2}2\sum_{\ell\geq0}\frac{\He_{2\ell}(0)}{(2\ell+1)!}\He_{2\ell+1}\left(\frac{x}{\s\sqrt n}\right)\frac{e^{-\frac{x^2}{2\s^2}}}{2\pi \sqrt{ n\s^2}}+\oo{\sqrt n}\nn\\
&=&\frac{d}{\s^2}\sum_{\ell\geq0}\left[2\int_0^\infty \frac{dy}{\sqrt{2\pi}} \frac{\He_{2\ell+1}(y)}{(2\ell+1)!}e^{-\frac{y^2}{2}}\right]\He_{2\ell+1}\left(\frac{x}\s\right)
\frac{e^{-\frac{x^2}{2\s^2}}}{ \sqrt{2\pi n\s^2}}+\oo{\sqrt n}\nn\\
&=&\frac{d}{\s^2}\sum_{\ell\geq0}\left[\int_0^\infty \frac{dy}{\sqrt{2\pi}} \sign(y) \frac{\He_{2\ell+1}(y)}{(2\ell+1)!}e^{-\frac{y^2}{2}}\right]\He_{2\ell+1}\left(\frac{x}{\s\sqrt n}\right)\frac{e^{-\frac{x^2}{2\s^2}}}{\sqrt{2\pi n\s^2}}+\oo{\sqrt n}\nn\\
&=&\frac{d}{\s^2}\sign (x)\frac{e^{-\frac{x^2}{2\s^2}}}{\sqrt{2\pi n\s^2}}+\oo{\sqrt n}\,.\nn
\eea
\end{proof}

This concludes the proof of the second lemma and therefore of the main theorem.
\end{proof}


\begin{thebibliography}{9}

\bibitem{C} C. Boldrighini, A. Marchesiello, C. Saffirio, {\em Weak dependence for a class of local functionals of Markov chains on $\Z^d$}, Methods Funct. Anal. Topology, {\bf 21}:302-314, (2015).

\bibitem{BP} C. Boldrighini, A. Pellegrinotti. {\em Random Walk on $\mathbb{Z}$ with One-Point Inhomogeneity},
Markov Proc. and Rel. Fields, \textbf{18}, 421-440, (2012).

\bibitem{bat} G. Szeg\"o. {\em Orthogonal Polynomials}, AMS Colloquium Publications, Vol. XXIII, Providence, Rhode Island (1933).

\bibitem{siberia} D. A. Korshunov. {\em Limit Theorems for General Markov Chains}, Siberian Math. Journal {\bf 42}:2 301-316, (2001). 

\bibitem{law} G. F. Lawler, V. Limic. {\em Random Walk: a Modern Introduction}, Cambridge University Press, Cambridge, (2010).

\bibitem{MZ} R. A. Minlos, E. A. Zhizhina. {\em Local Limit Theorem for a Non Homogeneous Random Walk on the Lattice}, Theory Probab. Appl. \textbf{39},  513-529, (1994).
 
\end{thebibliography}
\end{document}